\documentclass[]{article}
\newcommand{\Title}{Short-length routes in low-cost networks \emph{via}
  Poisson line patterns}

\usepackage[english]{babel}
\usepackage{calc, xspace}
\usepackage{amsmath, amsthm, amsfonts}
 \usepackage[mathcal,mathscr]{eucal}
\usepackage{url}

\usepackage[monochrome]{color}
\usepackage{graphicx}

 \newcommand{\href}[2]{#2}

\theoremstyle{plain}
\newtheorem{theorem}{Theorem}
\newtheorem{lem}{Lemma}[section]
\newtheorem{corollary}{Corollary}
\theoremstyle{definition}
\newtheorem{definition}{Definition}
\theoremstyle{remark}
\newtheorem{remark}{Remark}
\newtheorem{ex}{Example}[section]

 \newcommand{\Expect}[1]{\operatorname{\mathbb{E}}\left[#1\right]}
 \newcommand{\Indicator}[1]{\operatorname{\mathbb{I}}\left[#1\right]}
 
 \newcommand{\Law}[1]{\mathcal{L}\left({#1}\right)}
 \newcommand{\Prob}[1]{\operatorname{\mathbb{P}}\left[#1\right]}
 \newcommand{\Reals}{\mathbb{R}}

 \DeclareMathOperator{\Number}{\#}
 \DeclareMathOperator{\Steiner}{ST}

 \renewcommand{\d}{\text{\rm{d}}}
 
 \newcommand{\half}{\frac{1}{2}}
 \DeclareMathOperator*{\average}{average}
 \DeclareMathOperator{\config}{\textbf{x}}
 \DeclareMathOperator{\excess}{excess}
 \DeclareMathOperator{\len}{len}
 \DeclareMathOperator*{\ratio}{ratio}
 \DeclareMathOperator{\segment}{\underline{s}}

\newcommand{\thalf}{\tfrac{1}{2}}
\newcommand{\by}{\textbf{y}}
 
\begin{document}


 \title{\Title}
 \author{%
David J.~Aldous\thanks{Research supported by N.S.F Grant DMS-0203062}%
\; and Wilfrid S.~Kendall%
}
\date{\text{}}


 \maketitle

 \begin{abstract}
   In designing a network to link $n$ points in a square of area $n$,
   one might be guided by the following two desiderata.  First, the
   total network length should not be much greater than the length of
   the shortest network connecting all points.  Second, the average
   route length (taken over source-destination pairs) should not be much
   greater than the average straight-line distance.  How small can we
   make these two excesses?  Speaking loosely, for a non-degenerate
   configuration
 the total
   network length must be at least of order $n$ and the average straight-line
   distance must be at least of
 order $n^{1/2}$, so it seems implausible that
   a single network might exist
in which the excess over the first minimum is $o(n)$ and the excess over the
   second minimum is $o(n^{1/2})$.  But in fact one can do better:
   for an arbitrary configuration one can construct a network where
   the first excess is $o(n)$ and the second excess is almost
   as small as $O(\log n)$.  The construction is conceptually simple
   and uses stochastic methods:
   over the minimum-length connected network (Steiner tree)
   superimpose a sparse stationary and isotropic Poisson line
   process. Together with a few additions (required for technical reasons), the
   mean values of the excess for the resulting random 
   network satisfy the above asymptotics;
   hence a standard application of the probabilistic method guarantees
   the existence of deterministic networks as required
   (speaking constructively, such networks can be constructed using simple rejection sampling). The key
   ingredient is a new result about the Poisson line process.
   Consider two points at distance \(r\) apart, and delete from the
   line process all lines which separate these two points. The
   resulting pattern of lines partitions the plane into cells; the
   cell containing the two points has mean boundary length \(\approx
   2r + \text{constant}\times\log r\).  Turning to lower bounds, consider a sequence of networks in \([0, \sqrt{n}]^2\)
   satisfying a weak equidistribution assumption. We show that if the first
   excess is $O(n)$ then the second excess cannot be
   $o(\sqrt{\log n})$.
 \end{abstract}




\emph{Short title:} Lengths and costs in networks

 \section{Introduction}\label{sec:introduction}
 We start with a counter-intuitive observation and its motivation,
 which prompted us to probe more deeply into the underlying question.

 Consider \(n\) points (``cities'', say) in a square of area \(n\).
 Using the terminology of computer science, we are interested in both
 the \emph{worst-case} setting where the points are located
 arbitrarily in the square, and the \emph{average case} setting where
 the points are random, independent and uniformly distributed.
 Consider a connected network (a road network, say), made up of a
 finite number of straight line segments and linking these \(n\)
 points and perhaps other junction points.  Recall that the minimum
 length connected network on a configuration of points \(\config^n =
 \{x_1,\ldots,x_n\}\) is the \emph{Steiner tree}
 \(\Steiner(\config^n)\).

 It is well known and straightforward to prove
 \cite{Steele-1997,Yukich-1998} that in both the worst case and
 the average case the (mean) total network length
 \(\len(\Steiner(\config^n))\) grows as order \(O(n)\).  When
 designing a network, it is reasonable to regard total network length
 as a ``cost''.  A natural corresponding ``benefit'' would be the
 existence (in some average sense) of short routes between points.
 Let \(\ell(x_i,x_j)\) be the route-length (length of shortest path)
 between points \(x_i\) and \(x_j\) in a given network, and let
 \(|x_i-x_j|\) denote Euclidean distance (so
 \(\ell(x_i,x_j) \geq |x_i-x_j|\)). A good network should satisfy the
 following informal criterion:
 \begin{quote}
The \textbf{short routes property:}
   Averaging over pairs \((i,j)\) chosen uniformly at random, the route-length \(\ell(x_i, x_j)\) between
   points \(x_i\) and \(x_j\) is not
   much larger than the Euclidean distance \(|x_i-x_j|\).
 \end{quote}

 A first take on a statistic to measure this property for a connected
 network \(G(\config^n)\) is the \emph{ratio statistic},
based on averaging the
 ratios of network route-lengths \emph{versus} Euclidean distances. Consider
 a network \(G(\config^n)\) to be the configuration of points
 \(\config^n=\{x_1, \ldots, x_n\}\) together with a collection of line
 segments which combine to connect every \(x_i\) to
 every other \(x_j\).
 \begin{definition}[Ratio statistic]
\label{def:ratio-statistic}
Let \(\average_{(i,j)}\) denote the average over all distinct pairs
\((i,j)\). Then
   \begin{equation}
     \ratio(G(\config^n)) \quad=\quad
 \average_{(i,j)} \frac{\ell(x_i, x_j)}{|x_i-x_j|} - 1
     \quad\geq\quad 0\,.
     \label{eq:ratio-statistic}
   \end{equation}
 \end{definition}
%
%
     Consider a network \(G(\config^n)\) based on \(n\) uniform random
     points \(\config^n\subset[0,\sqrt{n}]^2\), having (say) twice the
     total length of the Steiner tree. Initially we speculated that in
     this case the expectation \(\Expect{\ratio(G(\config^n))}\)
     would at best converge to some strictly positive constant as
     \(n\to\infty\). However this intuition is wrong:
\begin{quote}
  \textbf{Counterintuitive observation} (see section
     \ref{sec:counterintuitive}). It is possible to construct
  networks over well-dispersed configurations whose total lengths  are
  greater than the corresponding Steiner tree lengths  by only an
  asymptotically negligible factor, but for which the ratio statistic
  converges to zero as total network length converges to infinity.
\end{quote}

These considerations were originally motivated by analysis of
real-world networks. Consider for example the ``core'' part of the UK
rail network; that part which links the \(40\) largest cities. Given a
statistic \(R\) designed to capture the ``short routes'' property, one
can then consider how closely the observed value of \(R\) approaches
optimality. Of course the real network has evolved according to a
complex historical process heavily influenced by topography;
nevertheless it is of interest to consider whether its value of \(R\)
is close to the minimum possible value of \(R\) taken over all possible
networks connecting the \(40\) cities but of no greater total length.

One is then led to ask
 what statistic \(R\)
might best capture the imprecisely expressed ``short routes''
property, and our consideration of \(n\) cities in an
idealised square \([0,\sqrt{n}]^2\) is designed to illuminate this
question. The above counterintuitive observation can be interpreted as
implying that the ratio statistic of Definition \ref{def:ratio-statistic}
is probably \emph{not} a good choice of
statistic, because we prove this observation by constructing networks
which are approximately optimal by this criterion and yet are plainly
rather different from many plausible real-world networks.
What
\emph{is} a good choice of statistic will be discussed in a companion
paper, along with some real-world examples.


Informally, the counter-intuitive observation suggests that we can
construct networks for configurations of \(n\) points which have total
network length exceeding that of the Steiner tree by just \(o(n)\),
and such that the average excess of network distance over Euclidean
distance is \(o(n^{1/2})\) (bearing in mind that average Euclidean distance
for ``evenly spread out'' configurations should be \(O(n^{1/2})\)). In
fact much more is true:
whatever the configuration of \(n\) points in \([0, \sqrt{n}]^2\) (hence, even in ``worst case'' scenarios)
we can construct such networks with average excess of network
distance over Euclidean distance barely more than \(O(\log n)\). This
we can work on an additive rather than a multiplicative scale:
\begin{definition}[Excess average length for a network]
\label{def:excess-route-length}
  The \emph{excess route length} for a network \(G(\config^n)\) is
  \begin{equation}
    \label{eq:excess-route-length}
    \excess\left(G(\config^n)\right)\quad=\quad
\average_{(i,j)}\left(\ell(x_i,x_j)-|x_i-x_j|\right)\,.
  \end{equation}
\end{definition}
\begin{theorem}[Upper bound on minimum excess network length]
\label{thm:upper-bound}
  For each \(n\) let \(\config^n\) be an arbitrary configuration of
  \(n\) points in a square of area \(n\). The following asymptotics
  hold for large \(n\):
  \begin{itemize}
  \item[(a)]
    \label{thm:upper-bound-A}
    Let \(w_n\to\infty\).  There exist networks \(G(\config^n)\)
    connecting up the points such that
    \begin{enumerate}
    \item[(i)] \(\len(G(\config^n))-\len(\Steiner(\config^n))=o(n)\);
    \item[(ii)] \(\excess(G(\config^n))=o(w_n\log n)\).
    \end{enumerate}
  \item[(b)]
    \label{thm:upper-bound-B}
    Let \(\varepsilon>0\).  There exist networks \(G(\config^n)\)
    connecting up the points such that
    \begin{enumerate}
    \item[(i)]
      \(\len(G(\config^n))-\len(\Steiner(\config^n))\leq\varepsilon
      n\);
    \item[(ii)] \(\excess(G(\config^n))=O(\log n)\).
    \end{enumerate}
  \end{itemize}
\end{theorem}
This result is proved in Sections \ref{sec:plp} and \ref{sec:upb}. The
idea is to build a hierarchical network. Details are given at the
start of Section \ref{sec:upb}, but here is a sketch.
At small scales routes use the underlying Steiner tree. At large
scales, routes use a sparse collection of randomly oriented lines (a
realisation of a stationary and isotropic \emph{Poisson line
  process}); this is the key ingredient that permits an excess of at
most \(o(w_n\log(n))\), respectively \(O(\log(n))\) (Section \ref{sec:plp}).  We believe that only
these two scales are needed, but to simplify matters (so as to avoid
non-elementary analysis of Steiner trees and geodesics in Poisson line
networks) we introduce an intermediate
scale consisting of a widely-spaced grid.  Thus a route from an originating city
navigates through the Steiner tree to a grid line and then along the
grid line to a line of the Poisson line process, and then navigates in
the reverse sense down to the destination city.  (For technical
reasons the discussion in Section
\ref{sec:upb} also introduces occasional small rectangles to permit
circumnavigation around Steiner tree ``hot-spots''). The key ingredient in the analysis is a calculation
concerning the Poisson line process,
%
%
which has separate interest as a result in stochastic geometry (Theorem
\ref{thm:asymptotic-mean-cell-boundary} below).  Consider
two points at distance \(r\) apart, and delete all lines from the line process
which separate these two points. The resulting pattern of lines
partitions the plane into cells; the cell containing the two points
has mean boundary length which for large \(r\) is asymptotic to \(2r +
\text{constant}\times\log r\).
%
%

Note that randomness arises only through use of the Poisson line
process to supply a relatively small number of long straight
connections; the point pattern \(\config^n\) is arbitrary. The
probabilistic method may now be used to prove the existence of a
non-random networks satisfying the asymptotics described in Theorem
\ref{thm:upper-bound}, based on applying Markov's inequality to the
expectations
\(\Expect{\len(G(\config^n))}-\len(\Steiner(\config^n))=o(n)\),
\emph{et cetera}.

For \emph{lower bounds} it is necessary to impose some condition on
the empirical distribution of the points in \(\config^n\), since
if all the points concentrate on a line then the excess is zero!
We need a quantitative condition on equidistribution of
points over a region, formalised via the following truncated
\emph{Vasershtein coupling} scheme.
\begin{definition}[Quantitative equidistribution condition]
\label{def:equidistribution}
Let \(\config^n\) for varying \(n\) form a sequence of configurations
in the plane, let \(\mu^n\) be a probability measure on the plane, and
and let \(L_n > 0\).  Say
\emph{
\(\config^n\) is \(L_n\)-equidistributed as \(\mu^n\)} if there exists a coupling of random variables
\((X_n,Y_n)\) such that
\begin{itemize}
\item[(a)] \(X_n\) has uniform distribution on the finite point-set
  \(\config^n\),
\item[(b)] \(Y_n\) has distribution \(\mu^n\),
\item[(c)] \(\Expect{ \min\left(1, \frac{|X_n - Y_n|}{L_n}\right)} \to 0\)
  as \(n \to \infty\).
\end{itemize}
\end{definition}
A sufficient condition for the following result is that \(\config^n\) is \(L_n\)-equidistributed as the uniform distribution on the square of area $n$, for some \(L_n = o(\sqrt{\log n})\).
The purpose of introducing the \emph{non}-uniform distribution
\(\mu^{n}\) in Definition \ref{def:equidistribution} is to permit us
to express Theorem
\ref{thm:lower-bound} below in terms of weaker and more local
conditions: for example a consequence of Theorem
\ref{thm:lower-bound}(b) is that we may replace the
\emph{uniform} reference distribution by any distribution \(\mu\) on
\( [0,1]^2\) with a continuous density component, rescaled to produce
a distribution \(\mu^{n}\) on \([0,n^{1/2}]^2\). In particular the
geometry of \([0,n^{1/2}]^2\) plays no r{\^o}le in this result.

We choose to express Definition \ref{def:equidistribution} in
stochastic terms purely for convenience of exposition. For example,
arguments using the connection of total variation to coupling show
that \(\config^n\) is \(L_n\)-equidistributed as the uniform
distribution on \([o,\sqrt{n}]^2\) if the following non-stochastic
condition is satisfied: for some sequence of numbers
\(\lambda_n\to\infty\) with \(\lambda_n/L_n\to0\) and \(n/\lambda_n^2\)
being integral,
\[
\frac{1}{n}\sum\left|\Number(\config^n\cap\text{box})-\lambda_n^2\right|
 \quad\to\quad0\,,
\]
with the sum being taken over \(n/\lambda_n^2\) boxes partitioning
\([o,\sqrt{n}]^2\) into cells of sidelength \(\lambda_n\). Thus a wide
range of possible point patterns can be seen to be
\(L_n\)-equidistributed in the above sense.

\begin{theorem}[Lower bound on minimum excess network length]
\label{thm:lower-bound}
Let \(\config^n\) be a configuration of \(n\) points in a square
\([0,\sqrt{n}]^2\).  Let \(L_n = o(\sqrt{\log n})\).
Suppose either
\begin{itemize}
\item[(a)] \(\config^n\) is \(L_n\)-equidistributed as the uniform distribution on the square of area $n$;
\item[] or (more generally)
\item[(b)] for some fixed \(\rho\) and \(\varepsilon\), there is a
  subcollection \(\by^{k(n)}\) of \(k(n)\) points, all lying in a
  disk \(D_n\) of area \(\pi\rho n\), such that \(k(n)>\pi\rho
  n\varepsilon\), and such that \(\by^{k(n)}\) is
  \(L_n\)-equidistributed as the uniform distribution on \(D_n\).
\end{itemize}
Let \(G(\config^n)\) be a network based on the full collection of
\(n\) points.  If
\(\len(G(\config^n))/n\) remains bounded as \(n\to\infty\), then
\begin{equation}
  \label{eq:lower-bound}
  \excess(G(\config^n))\quad=\quad \Omega(\sqrt{\log n})\,.
\end{equation}
(Thus, \(\liminf_{n\to\infty}\excess(G(\config^n))/\sqrt{\log n}>0\).)
\end{theorem}

Configurations \(\config^n\) produced by independent uniform sampling
from \([0,\sqrt{n}]^2\) satisfy the conditions of this theorem (see
Remark \ref{rem:reassuring-exercise}), but so will many other
configurations exhibiting both clustering and repulsion.
The proof of the theorem is given in Section \ref{sec:lwb}, and exploits a tension
between the two following facts:
\begin{enumerate}
\item[(a)]
\label{item:parallel}
An efficient route between \(x_i\) and \(x_j\) must run approximately
parallel to the Euclidean geodesic, and hence will tend to make
almost orthogonal intersections with random segments perpendicular to
this geodesic.
\item[(b)]
\label{item:not-parallel}
On the other hand, the equidistribution condition means that two
points \(x_i\) and \(x_j\) randomly chosen from the subcollection must
be nearly independent uniform draws from \(D_n\), which permits the
derivation of \emph{upper bounds} on the probability of nearly
orthogonal intersections of the form given in fact (a).
\end{enumerate}

Finally, one might hope to improve the result by imposing a more restrictive
assumption than the requirement that \(\len(G(\config^n)/n\) remains bounded as
\(n\to\infty\).  This requirement is weaker than either of our two alternative
assumptions on \(\len(G(\config^n))-\len(\Steiner(\config^n))\) in the
upper bound (since \(\len(\Steiner(\config^n))=O(n)\)). However we are
unable to improve (\ref{eq:lower-bound}) under either of the two stronger
assumptions.

 \section{The Poisson line process network}\label{sec:plp}
 Our upper bound on minimal \(\excess\left(G(\config^n)\right)\) is
 based on a result from stochastic geometry (Theorem
 \ref{thm:asymptotic-mean-cell-boundary} below) which is of independent
 interest.

\begin{figure}[ht]
     \centering
     \includegraphics[width=3in]{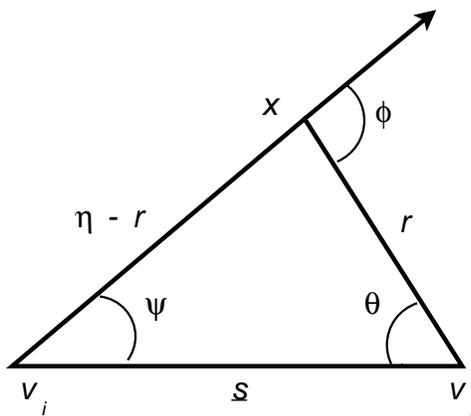}
     \caption{Definition of \(\eta\) and
       \(\phi\). Note that \(\phi\) is the sum of the two interior
       angles \(\psi\) and \(\theta\).}
\label{fig:triangle}
   \end{figure}

 Recall that a Poisson line process in the plane \(\Reals^2\) is
 constructed as a Poisson point process whose points lie in the space
 which parametrises the set of lines in the plane.
 We will consider
 only undirected lines, which will be
 parametrised by \((r,\theta)\in\Reals\times[0,\pi)\)
 where \(r\) is the
 signed distance from the line to a reference point and
 \(\theta\) is the angle the line makes with a reference axis. A
stationary and
 isotropic Poisson line process has intensity measure invariant under
 rotations and translations of \(\Reals^2\):
 a stationary and isotropic Poisson line process \(\Pi\) of unit
 intensity is one for which the number of lines of \(\Pi\) hitting a
 unit segment has expectation \(1\) (further facts about Poisson line
 processes may be found in \cite[Chapter
 8]{StoyanKendallMecke-1987}). We are interested in the cell
 containing two fixed points which is formed by
 the lines of \(\Pi\) that do not separate the two points, because this
 can be used as the efficient long-distance part of a network route between
 the two points (see Lemma \ref{med-large}).
Theorem
 \ref{thm:asymptotic-mean-cell-boundary} establishes an asymptotic
upper bound for
 the length of the mean cell perimeter in case of wide separation
 between the two points; we prepare for this by using a Buffon
 argument to derive an exact double-integral expression for the mean
 cell perimeter length:
   \begin{theorem}[Mean perimeter length]\label{thm:mean-cell-boundary}
     Let \(\Pi\) be a stationary and isotropic Poisson line process of
     unit intensity.    Fix two points \(v_i\), \(v_j\) which are distance \(m\) apart.  Delete the lines of \(\Pi\) which separate the
     two points \(v_i\), \(v_j\). The remaining line pattern partitions the plane:
     the cell \(\mathcal{C}(v_i,v_j)\) containing the two fixed points has
     mean perimeter \(\Expect{ \len \partial\mathcal{C}(v_i,v_j)}= 2m
     + J_m\),
where \(J_m\) is given by
     the double integral
     \begin{multline}
       \label{eq:double-integral}
       J_m \quad=\quad
       \Expect{\len\partial\mathcal{C}(v_i,v_j)}-2m 
\quad=\quad
       \half\iint\limits_{\Reals^2}
       \left(\phi-\sin\phi\right)
       \exp\left(-\thalf(\eta-m)\right)
       {\d}x\,.
     \end{multline}
     Here \(\eta=\eta(x)\) is a sum of distances
     \(|v_i-x|+|v_j-x|\), while \(\phi=\phi(x)\) is the
     exterior angle at \(x\) of the triangle with vertices \(x\),
     \(v_i\), \(v_j\) (see Figure \ref{fig:triangle}).
   \end{theorem}
\begin{proof}
 This proof can be phrased in terms of measure-theoretic stochastic
 geometry, using the language of Palm distributions and Campbell
 measure. Since we deal only with constructions based on Poisson
  processes, we are able to adopt a less formal but more transparent
  exposition, for the sake of a wider readership.

  Let \(\segment\) be the line segment of length \(m\) with end-points
  \(v_i\), \(v_j\). The idea of the proof is to measure
  \(\Expect{\len\partial\mathcal{C}(v_i,v_j)}\) by computing the expected
  number of hits on \(\partial\mathcal{C}(v_i,v_j)\)
made by an \emph{independent}
  homogeneous isotropic Poisson line process \(\widetilde{\Pi}\),
  again of unit intensity. Each hit corresponds to one of the points in the
  \emph{intersection point process}
  \(\mathcal{X}=\{\iota(\ell,\widetilde{\ell}):\ell\in\Pi,
  \widetilde{\ell}\in\widetilde{\Pi}\}\), where
  \begin{equation}
    \label{eq:intersection-map}
    \iota(\ell, \widetilde{\ell})\quad=\quad
    \begin{cases}
      x & \text{ if }\ell\cap\widetilde{\ell}=\{x\}\,,\\
      \text{undefined} & \text{ if } \ell, \widetilde{\ell}
\text{ are parallel.}
    \end{cases}
  \end{equation}
  Note that with probability \(1\) the intersection point \(\iota(\ell,
  \widetilde{\ell})\) is defined for all \(\ell\in\Pi\),
  \(\widetilde{\ell}\in\widetilde{\Pi}\).

Not all intersection points \(x\in\mathcal{X}\) correspond to hits on
\(\partial\mathcal{C}(v_i,v_j)\). The condition for \(x=\iota(\ell, \widetilde{\ell})\in\mathcal{X}\) to represent
a hit on
\(\partial\mathcal{C}(v_i,v_j)\) is that \(\ell\) should not hit
\(\segment\) (for otherwise it cannot be involved in the construction
of \(\partial\mathcal{C}(v_i,v_j)\)) and that 
\(x\) is not separated from
\(\segment\) by any line from \(\Pi\setminus\{\ell\}\). Recall that the Slivynak
theorem \cite[\S4.4, example 4.3]{StoyanKendallMecke-1987} implies
that \(\Pi\setminus\{\ell\}\) conditional on \(\ell\in\Pi\) is itself
a homogeneous isotropic unit-rate Poisson line process. Consequently,
under the condition that
\({\ell}\) does not hit \(\segment\),
the conditional probability of
 \(x=\iota(\ell, \widetilde{\ell})\in\mathcal{X}\)
representing a hit on \(\partial\mathcal{C}(v_i,v_j)\) is equal to the probability \(p(x)\) of
there being no line in \(\Pi\) which cuts both the segment from \(v_i\) to
\(x\) and the segment from \(v_j\) to \(x\).

A classic counting argument from stochastic geometry then reveals that
\begin{equation}
\label{eq:stochastic-geometry-counting}
  p(x)\;=\;
  \exp\left(-\thalf\left(|v_i-x|+|v_j-x|-m\right)\right)
\quad=\quad\exp\left(-\thalf(\eta-m)\right)\,.
\end{equation}
Accordingly, if \(\nu\) is the intensity of the point process
\(\mathcal{X}\) then we may compute the mean number of hits on
\(\partial\mathcal{C}(v_i,v_j)\) as
\begin{multline}
\iint_{\Reals^2}\nu\Prob{\ell\not\Uparrow\segment |
x=\iota(\ell,\widetilde{\ell})\in\mathcal{X}}\exp\left(-\thalf(\eta-m)\right)
{\d}x\\
\quad=\quad
2m+
  \iint_{\Reals^2}
\nu
\Prob{\ell\not\Uparrow\segment, \widetilde{\ell}\not\Uparrow\segment |
x=\iota(\ell,\widetilde{\ell})\in\mathcal{X}}\exp\left(-\thalf(\eta-m)\right)
{\d}x .
\end{multline}
Here ``\(\ell\not\Uparrow\segment\)'' stands for ``the line \(\ell\)
does not hit \(\segment\)'' -- noting that the conditioning in this
context forces
the Poisson line \(\ell\) to pass through \(x\) but does not fix its
orientation -- and on the right-hand side the summand \(2m\) corresponds to the fact that
hits of \(\widetilde{\Pi}\) on \(\segment\) count as automatic double hits on
\(\partial\mathcal{C}(v_i,v_j)\).

Condition on \(x=\iota(\ell,\widetilde{\ell})\in\mathcal{X}\) (which
is to say, condition on there being Poisson lines \(\ell\in\Pi\),
\(\widetilde{\ell}\in\widetilde{\Pi}\) both passing through \(x\)) and
consider
\begin{itemize}
\item[(a)] the angle \(\xi_1\) made by \(\ell\) with the line through
  \(v_i\) and \(x\);
\item[(b)] the angle \(\xi_2\) between \(\ell\) and \(\widetilde{\ell}\).
\end{itemize}
By isotropy of \(\Pi\) the random angle \(\xi_1\) is
Uniform(\(0,\pi)\). Conditional on \(\xi_1\) and more generally on
\(\Pi\) with an \(\ell\in\Pi\) passing through \(x\), the intersection
of \(\widetilde{\Pi}\) with \(\ell\) is a Poisson point process on
\(\ell\) of unit intensity. Moreover if the intersection points are
marked with angles of intersection \(\xi_2\) then the mark \(\xi_2\)
has mark density \(\tfrac{1}{2}\sin\xi_2\) over \(\xi_2\in[0,\pi)\) (consider
the length of the silhouette of a portion of \(\ell\) viewed at angle
\(\xi_2\)). Hence the conditional distribution of \(\xi_2\) for
\(x=\iota(\ell,\widetilde{\ell})\) has
density \(\tfrac{1}{2}\sin\xi_2\) over \(\xi_2\in[0,\pi)\), and so
we can compute (working with \(\xi_2\) modulo \(\pi\))
\begin{multline}
  \Prob{\ell\not\Uparrow\segment, \widetilde{\ell}\not\Uparrow\segment |
x=\iota(\ell,\widetilde{\ell})}\\
\quad=\quad
\frac{1}{\pi}\int_0^{\phi}\left(\int_{-\xi_1}^{\phi-\xi_1}\frac{|\sin\xi_2|}{2}\d\xi_2\right)\d\xi_1
\quad=\quad \frac{\phi-\sin(\phi)}{\pi}
\end{multline}
where 
\(\phi=\theta+\psi\) is the exterior angle at \(x\) of the triangle formed by \(x\), \(v_i\), \(v_j\) (see Figure
\ref{fig:triangle}).

Finally the intensity \(\nu\) of \(\mathcal{X}\) can be computed as
\(\tfrac{\pi}{2}\), for example by computing the mean number of hits
of the unit disk by \(\Pi\), then by computing the average length of
the intersection of the disk with a line of \(\Pi\) conditional on
that line hitting the disk. Thus
\begin{multline}
  J_m\quad=\quad\Expect{\len(\partial\mathcal{C}(v_i,v_j))}-2m\\
\quad=\quad
\nu\iint_{\Reals^2}
\Prob{\ell\not\Uparrow\segment, \widetilde{\ell}\not\Uparrow\segment |
x=\iota(\ell,\widetilde{\ell})\in\mathcal{X}}\exp\left(-\thalf(\eta-m)\right)
{\d}x\\
\quad=\quad
\frac{1}{2}\iint_{\Reals^2}
\left(\phi-\sin\phi\right)\exp\left(-\thalf(\eta-m)\right)
{\d}x
\end{multline}
as required.
\end{proof}

We now state and prove the main result of this section: an \(O(\log m)\)
upper bound on the mean perimeter excess length \(J_m\).
   \begin{theorem}[Asymptotic upper bound on mean perimeter length]
\label{thm:asymptotic-mean-cell-boundary}
The mean perimeter excess length \(J_m\) is subject
to the following asymptotic upper bound:
     \begin{equation}
       \label{eq:leading-asymptotics}
       J_m \quad\leq\quad O(\log m) \qquad\text{as }m\to\infty\,.
     \end{equation}
   \end{theorem}
\begin{proof}
Without loss of generality,
place the points $v_i$ and $v_j$ at \((-\tfrac{m}{2},0)\) and \((\tfrac{m}{2},0)\).
  The double integral in \eqref{eq:double-integral} possesses mirror
  symmetry about each of the two axes, so we can write
  \begin{multline}
    \label{eq:mirror2}
    J_m\quad=\quad
2\iint\limits_{[0,\infty)^2}
       \left(\phi-\sin\phi\right)
       \exp\left(-\thalf(\eta-m)\right)
       {\d}x\\
\quad=\quad
2\int_0^{\pi/2}\int_0^{\frac{m}{2}\sec\theta}
\left(\phi-\sin\phi\right)\exp\left(-\thalf(\eta-m)\right)r\,\d r\,\d \theta+\\
+
2\int_{\pi/2}^{\pi}\int_0^\infty
\left(\phi-\sin\phi\right)\exp\left(-\thalf(\eta-m)\right)r\,\d r\,\d\theta
  \end{multline}
  (using polar coordinates \((r,\theta)\) about the second point
  \(v_j\) located at \((\tfrac{m}{2},0)\)).  The integrand in the
  second summand is dominated by
  \(\pi\exp\left(-\tfrac{r}{2}\right)r\), which is integrable over
  \((r,\theta)\in(0,\infty)\times(\tfrac{\pi}{2},\pi)\). (In this
  region geometry shows that \(\eta-m>r(1-\cos\theta)\geq r\).) Thus
  we can apply Lebesgue's dominated convergence theorem to deduce that
  the second summand is \(O(1)\) as \(m\to\infty\), hence may be
  neglected.

In fact we can also show that part of the first summand generates an
\(O(1)\) term: the dominated convergence theorem can be applied for
any \(\varepsilon\in(0,\pi/2]\) to show that
\[
2\int_0^{\pi/2}\int_\varepsilon^{\frac{m}{2}\sec\theta}
\left(\phi-\sin\phi\right)\exp\left(-\thalf(\eta-m)\right)r\,\d r\,\d\theta
\quad=\quad O(1)\,,
\]
since the integrand is dominated by
\(\pi\exp\left(-\tfrac{r}{2}(1-\cos\theta)\right)r\) over the region
\((r,\theta)\in(0,\infty)\times(\varepsilon,\tfrac{\pi}{2})\) (in this region
geometry shows that
\(\eta-m>r(1-\cos\theta)>r(1-\cos\varepsilon)\)). Thus
for fixed \(\varepsilon\in(0,\tfrac{\pi}{2})\)
as \(m\to\infty\) we have the asymptotic expression
\[
J_m\quad=\quad
2\int_0^{\varepsilon}\int_0^{\frac{m}{2}\sec\theta}
\left(\phi-\sin\phi\right)\exp\left(-\thalf(\eta-m)\right)r\,\d r\,\d\theta
+O(1)\,,
\]
where the implicit constant of the \(O(1)\) term depends on the choice
of\(\varepsilon>0\).

Now in the region where
\(0<\theta<\varepsilon\) and \(0<r<\tfrac{m}{2}\sec\theta\)
 we
know that \(\phi<2\theta<2\varepsilon\), and moreover \(\phi-\sin\phi\) is
an increasing function of \(\phi\).
 Therefore there is
a constant \(C_\varepsilon\), converging to zero as
\(\varepsilon\to0\),
 such that in this region
\begin{equation*}
  \phi-\sin\phi \quad\leq\quad 2\theta - \sin(2\theta) \quad\leq\quad
  \frac{C_\varepsilon}{8} \  \frac{(2 \theta)^3}{6}
\quad\leq\quad C_\varepsilon \frac{1-\cos\theta}{3}\sin\theta\,.
\end{equation*}
Hence (as \(m\to\infty\) for fixed \(\varepsilon>0\))
\begin{align*}
  &2\int_0^{\varepsilon}\int_0^{\frac{m}{2}\sec\theta}
\left(\phi-\sin\phi\right)\exp\left(-\thalf(\eta-m)\right)r\d r\d\theta\\
&\quad\leq\quad
\tfrac{2}{3} C_\varepsilon\int_0^{\varepsilon}\int_0^{\frac{m}{2}\sec\theta}
\left(1-\cos\theta\right)\sin\theta
\exp\left(-\tfrac{r}{2}(1-\cos\theta)\right)r\d r\d\theta\\
&\quad=\quad
\tfrac{8}{3} C_\varepsilon\int_0^{\varepsilon}
\left(\int_0^{\frac{m}{4}(\sec\theta-1)}
s e^{-s}
\d s\right)
\frac{\sin\theta\,\d\theta}{1-\cos\theta}
\qquad(\text{using }s=\tfrac{r}{2}(1-\cos\theta))\\
&\quad\leq\quad
\tfrac{8}{3} C_\varepsilon\int_0^{\frac{m}{4}(\sec\varepsilon-1)}
\left(\int_0^v
s e^{-s}
\d s\right)
\frac{1}{1+4v/m}\,\frac{\d v}{v}
\qquad(\text{using }v=\tfrac{m}{4}(\sec\theta-1))\\
&\quad\leq\quad
\tfrac{8}{3} C_\varepsilon \log\left(\tfrac{m}{4}(\sec\varepsilon-1)\right)
+O(1)\,.
\end{align*}
\end{proof}

\begin{figure}[ht]
  \centering
  \includegraphics[width=3in]{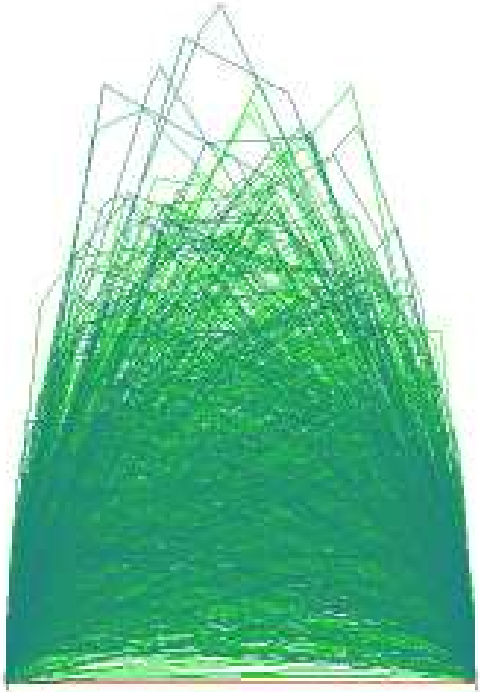}
  \caption{Simulation of semi-perimeters for \(1000\) independent
    cells for unit-rate Poisson line process, with points located
    at distance \(10^8\) units apart. The figure is subject to
    vertical exaggeration: \(y\)-axis is scaled at \(10^4\) times
    \(x\)-axis. Empirical mean excess semi-perimeter is \(27.63\) with
    standard error \(\pm0.28\), \emph{versus} predicted mean excess
    semi-perimeter \(27.5528\) (using \(o(1)\)-asymptotics).}
\label{fig:100n1000i}
\end{figure}

\begin{remark}
\label{rem:further-asymptotic-mean-cell-boundary}
More careful analysis yields useful \(o(1)\)-asymptotics: in
fact it can be shown that as \(m\to\infty\) so
     \begin{equation}
       \label{eq:further-asymptotic-mean-cell-boundary}
  J_m \quad=\quad \tfrac{8}{3}\left(\log m+\gamma+\tfrac{5}{3}\right)+o(1)\,,
     \end{equation}
     where \(\gamma\) is the Euler-Mascheroni constant:
     \(\gamma=\lim_{m\to\infty}\left(\left(\sum_1^m
         \tfrac{1}{r}\right)-\log m\right)\).  These
     \(o(1)\)-asymptotics show very good agreement with simulation:
     see for example the simulation reported in the legend of Figure
     \ref{fig:100n1000i}.
\end{remark}
 \section{A low-cost network with short routes}\label{sec:upb}
\begin{figure}[t]
  \centering
  \includegraphics[width=3in]{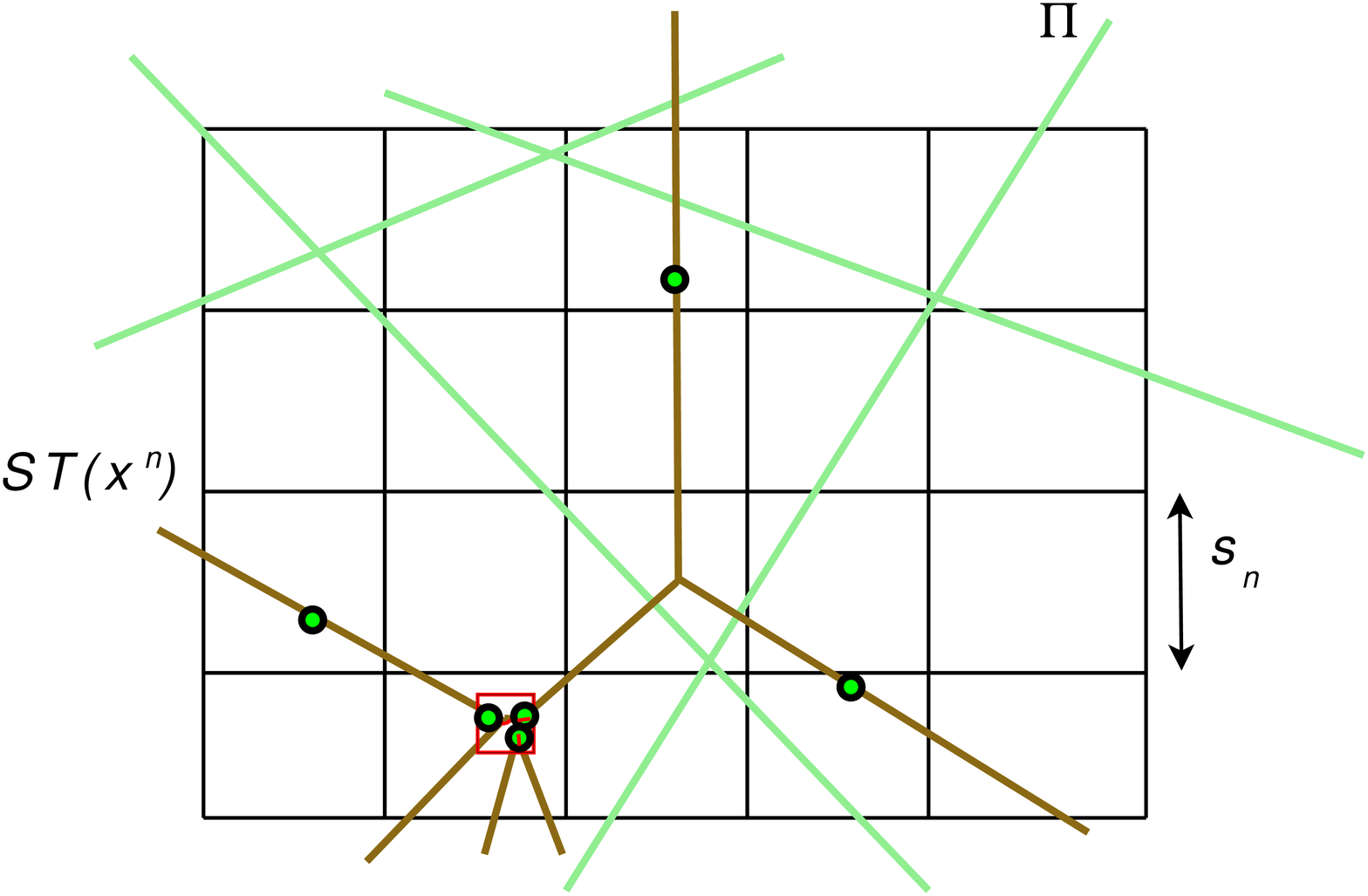}
  \caption{Illustration of construction of network to deliver an upper
    bound on mean excess route-length. Points are indicated by
    small circles. In this figure there is just one hot-spot cell.}
\label{fig:upb-1}
\end{figure}

 In this section we prove Theorem \ref{thm:upper-bound}: for a given
 configuration \(\config^n\subset[0,\sqrt{n}]^2\) we construct
 networks \(G(x^n)\) for which both
 \(\len(G(\config^n))-\len(\Steiner(\config^n))\) and
 \(\excess(G(\config^n)\) are small. The network is constructed by
 augmenting the Steiner tree network \(\Steiner(\config^n)\) in a
 hierarchical manner. The construction is stochastic: we construct a
 random augmentation for which the mean values of these excess values obey
 the desired asymptotics and then apply the probabilistic method to
 establish existence of the desired non-stochastic networks. Working
 from the largest scale downwards, we construct
\begin{enumerate}
\item a \emph{stationary and isotropic Poisson line process \(\Pi\) of
    intensity \(\eta\)}, where \(\eta\) will be small: note that this
  can be constructed from a unit intensity process by scaling by a
  magnification factor of \(1/\eta\). A
  simple computation \cite[Section 8.4]{StoyanKendallMecke-1987}
  shows that the mean total length of the intersection of the
  resulting line pattern with \([0,\sqrt{n}]^2\) equals \(\pi\eta n\).
\item A \emph{medium-scale rectangular grid with cell side-length \(s_n \sim (\log
    n)^{1/3}\).} Total length of this grid in \([0,\sqrt{n}]^2\) is
  bounded above by
\[
2 (1 + \tfrac{\sqrt{n}}{s_n}) \sqrt{n}
\quad=\quad o(n)\,.
\]
\item The \emph{Steiner tree \(\Steiner(\config^n)\)}.
\item A small number (at most \(n/2\)) of \emph{small hot-spot cells}
  based on a small-scale rectangular grid with cell side-length
  \(t_n \sim \tfrac{1}{(\log n)^{1/6}}\). A cell in this grid is
  described as
  a \emph{hot-spot cell} if it contains two or more points.
These hot-spot cells are used
  to by-pass regions where the Steiner tree might become complicated and
  expensive in terms of network traversal. We add further small segments
  connecting each hot-spot cell perimeter to points within the
  hot-spot cell. Total length of these additions can be bounded by
\[
4 \frac{n}{2} t_n + n \frac{t_n}{2}
\quad=\quad o(n)\,.
\]
\end{enumerate}
Thus the mean excess length of this augmented network is
\(o(n)+\pi\eta n\). The construction is illustrated in Figure
\ref{fig:upb-1}.  Note that we can choose $s_n$ and $t_n$ such that
$n^{1/2}/s_n$ and $s_n/t_n$ are integers, so that the small-scale
lattice is a refinement of the medium-scale lattice, which itself
refines the square \([0,\sqrt{n}]^2\).

\subsection{Worst-case results for Steiner trees}
\label{sec:steiner-worst-case}
We first record two elementary results on Steiner trees. The first
result bounds the length of a Steiner tree in terms of the square-root
of the number of points (for the planar case).
\begin{lem}\label{lem:steiner-1}
  Consider a configuration \(\config^k\) of \(k\) points in a square
  of side \(r\): there is a constant \(C_1\) not depending on \(k\) or
  \(r\) such that
\begin{equation}
  \len\left(\Steiner(\config^k)\right)\quad\leq\quad C_1 \sqrt{k} r\,.
\label{eq:steiner-1}
\end{equation}
\end{lem}
\begin{proof}
  See \cite[Section 2.2]{Steele-1997}.
\end{proof}

The second result provides a local bound on length contributed by a larger
Steiner tree in a small square containing a fixed number of points.
\begin{lem}\label{lem:steiner-2}
  Consider the Steiner tree \(\Steiner\left(\config^n\right)\) for an
  \emph{arbitrary} configuration \(\config^n\) in the plane. Let \(G\)
  be the restriction of the network \(\Steiner\left(\config^n\right)\)
  to a fixed open square of side-length \(t\). Suppose \(k\)
  points \(x_1\), \ldots, \(x_k\) of the configuration \(\config^n\)
  lie within the square. Then
  \begin{equation}
    \label{eq:steiner-2}
    \len(G)\quad\leq\quad t \left(4 + C_1\sqrt{k+1}\right)\,.
  \end{equation}
\end{lem}
\begin{proof}
  Let \(y_1\), \ldots, \(y_m\) be the locations at which
  \(\Steiner\left(\config^n\right)\) crosses into the interior of the
  square. (Note: \(m=0\) is possible if \(\{x_1, \ldots,
  x_k\}=\config^n\): in this case choose \(y_1\) arbitrarily from the
  perimeter of the square.) Then
  \begin{align*}
    \len(G)&\quad\leq\quad \len(\Steiner(\{x_1, \ldots,
  x_k, y_1, \ldots, y_m\}))
\qquad&\text{by minimality of \(\Steiner\left(\config^n\right)\),}\\
&\quad\leq\quad
\len(\Steiner(\{x_1, \ldots,
  x_k, y_1\})) +4t
\qquad&\text{using square perimeter,}\\
&\quad\leq\quad
t \left(4 + C_1\sqrt{k+1}\right)
\qquad&\text{using the previous lemma.}
  \end{align*}
\end{proof}

\subsection{Route-lengths in the medium-large network}
The part of the construction involving the medium-scale grid and the Poisson line process is useful in variant problems, so we separate out the following estimate involving these ingredients.
\begin{lem}\label{med-large}
Let $n^{1/2}/s_n$ be an integer.
Consider the superposition of the rectangular grid with cell side-length $s_n$ and the
Poisson line process of intensity $\eta$, intersected with the square $[0,n^{1/2}]^2$.
Let $v_i, v_j$ be vertices of the grid.
Then
\[ \Expect{ \mbox{route-length $v_i$ to $v_j$}}
\leq |v_i-v_j|
+ C_2  \tfrac{1}{\eta} \log (\eta \sqrt{2n}) \]
for an absolute constant $C_2$.
\end{lem}
\begin{proof}
Let
\(\mathcal{C}(v_i,v_j)\) be the cell of \(\Pi\) containing \(v_i\) and
\(v_j\) (having deleted lines from \(\Pi\) which separate
 \(v_i\) from \(v_j\)). Let
\(R(v_i,v_j)\) be the rectangle bounded by \(v_i\) and \(v_j\); then
by convexity the
route-length from $v_i$ to $v_j$
is bounded above by
\[
\half \len \partial\left(R(v_i,v_j)\cap\mathcal{C}(v_i,v_j)\right)
\quad\leq\quad\half \len \partial\mathcal{C}(v_i,v_j)\,,
\]
whose mean value can be computed by recognising that the Poisson line
process is a rescaled version of a homogeneous isotropic unit rate
Poisson line process. Hence by scaling the asymptotic upper bound of Theorem
\ref{thm:asymptotic-mean-cell-boundary} we have
\begin{multline*}
\Expect{\half \len \partial\left(R(v_i,v_j)\cap\mathcal{C}(v_i,v_j)\right)}
-|v_i-v_j|\\
\quad\leq\quad
O\left(\frac{1}{\eta}\log\left(\eta |v_i-v_j|\right)\right)\,
\quad=\quad
O\left(\frac{1}{\eta}\log\left(\eta \sqrt{2n}\right)\right) .
\end{multline*}

\end{proof}

\subsection{Navigating the augmented network}\label{sec:moving}
We now explain how to move from points of \(\config^n\) up to a vertex of the medium-scale grid.

Given \(x_i\in\config^n\), if this is in one of the hot-spot cells
then move to the perimeter of the hot-spot cell and thence to a
suitable point of departure on the perimeter, with route-length at most
\(\tfrac{5}{2}t_n\). Now move along the Steiner tree within the
relevant medium-scale grid box to the box perimeter; however by-pass
all hot-spot cells. There are \(\left(s_n/t_n\right)^2=\left((\log
  n)^{1/3}(\log n)^{1/6}\right)^2 =\log n\) small squares each of
which involves a route-length of either \(2t_n\) (if a hot-spot box which
will be by-passed) or \(t_n(4+C_1\sqrt{2})\) (if not, by Lemma
\ref{lem:steiner-2}). Hence the total trip to the medium-scale grid
box perimeter (including emergence from the initial hot-spot, if
required) has length at most
\[
\tfrac{5}{2}t_n+
t_n(4+C_1\sqrt{2}) \times {s_n^2}/{t_n^2}\quad \sim \quad
\tfrac{5}{2}t_n+(4+C_1\sqrt{2})\times(\log n)^{5/6}
\quad=\quad o(\log n)\,.
\]
Furthermore the route length from perimeter to vertex of medium-scale grid box is at
most
\(
\thalf s_n \sim  \tfrac{1}{2}(\log n)^{1/3} = o(\log n)\,.
\)
So for each $x_i$ there is a medium-scale grid vertex $v_i$ for which
\(
\mbox{route-length from $x_i$ to $v_i$ is } o(\log n)\) .
Combining with Lemma \ref{med-large} and noting that
the medium-scale grid geometry forces
\(|v_i-v_j|\leq|x_i-x_j|+2 \tfrac{s_n}{\sqrt{2}}\), we find
\[
\Expect{ \mbox{route-length from $x_i$ to $x_j$}} - |x_i-x_j| \leq
\sqrt{2}s_n + o(\log n) + C_2 \tfrac{1}{\eta}\log\left(\eta \sqrt{2n}\right) . \]
Averaging over the points of \(\config^n\), it follows that the
dominant contribution comes from the cell semi-perimeters, and indeed
\[
\Expect{\excess(G(\config^n))}\quad\leq\quad
O\left(\tfrac{1}{\eta}\log\left(\eta \sqrt{2n}\right)\right)\,,
\]
at a cost in terms of network length which exceeds
\(\len(\Steiner(\config^n))\) by a stochastic quantity of mean
\(\pi\eta n+o(n)\).

The two different results of Theorem \ref{thm:upper-bound} follow by
choosing \(\eta\) to behave in two different ways:
\begin{enumerate}
\item[(a)] either \(\eta\to0\), \(\eta w_n\to\infty\),
\item[(b)] or \(\eta=\varepsilon>0\).
\end{enumerate}
In either case we can apply the probabilistic method to exhibit
existence of the required deterministic networks for cases (a) and (b)
of Theorem \ref{thm:upper-bound}. For example in case (a) it is then
the case that
\(\Expect{\len(G(\config^n))-\len(\Steiner(\config^n))}\leq n c_n\)
and \(\Expect{\excess(G(\config^n))}\leq c_n w_n\log n\) for some
\(c_n\to0\). But then for any fixed \(n\) we can apply Markov's
inequality: \(\Prob{\len(G(\config^n))-\len(\Steiner(\config^n))> 3n
c_n}\leq\tfrac{1}{3}\) and \(\Prob{\excess(G(\config^n))>3 c_n w_n\log
n}\leq\tfrac{1}{3}\). Hence there is positive probability that the
random network satisfies both
\(\len(G(\config^n))-\len(\Steiner(\config^n))\leq 3n c_n\) and
\(\excess(G(\config^n))\leq3 c_n w_n\log n\), hence such a network
exists for each \(n\).

We can view these applications of Markov's inequality as indicating a
simple rejection sampling algorithm to be used to generate the
required sequence of networks.

 \section{A lower bound on average excess route-length}\label{sec:lwb}
 In this section we prove Theorem \ref{thm:lower-bound}. The proof is
 divided into four parts. Firstly (Subsection \ref{sec:lwb-reduction})
 we show how to reduce the problem to an analogous case in which the
 excess is computed for two random points drawn independently and
 uniformly from the whole disk \(D_n\) of area \(\pi\rho n\)
given in condition (b) of the
 theorem. Then (Subsection \ref{sec:lwb-geodesic-slope}) we show that
 the network geodesic must run almost parallel to the Euclidean
 geodesic if the excess is small. On the other hand (Subsection
 \ref{sec:lwb-uniform}) we can use the uniformity of the two random
 points to control the extent to which network segments can run
 both close to and nearly parallel to the Euclidean geodesic. Finally
 (Subsection \ref{sec:calculations}) we use the opposing estimates of
 Subsections \ref{sec:lwb-geodesic-slope} and \ref{sec:lwb-uniform} to
 derive a proof of the theorem using the method of contradiction.

 \subsection{Reduction to case of a pair of uniformly random points}
\label{sec:lwb-reduction}
First we indicate how condition (a) of Theorem \ref{thm:lower-bound}
implies condition (b). Under condition (a) we can use the coupling
between \(X_n\) and \(Y_n\) to show that \(\Number\{\config^n\cap
D_n\}/n\to\pi\rho\): therefore for large \(n\) the number of
points in \(\config^n\cap D_n\) is approximately \(\pi\rho n\). On the other hand
the same coupling can be used to bound the total variation distance
between the two conditional distributions \(\Law{Y_n | X_n\in D_n}\)
and \(\Law{Y_n | Y_n\in D_n}=\text{Uniform}(D_n)\), and to show that
this bound tends to zero. We can then use rejection sampling techniques
to couple \(\Law{Y_n | X_n\in D_n}\) and
\(\text{Uniform}(D_n)\) so that the truncated Vasershtein distance
tends to zero as \(n\to\infty\); as the distance is a metric we can combine this
coupling with the (conditioned) coupling of \(\Law{X_n|X_n\in D_n}\)
and \(\Law{Y_n|X_n\in D_n}\) to obtain a coupling which satisfies
condition \((b)\).

We now note that it is sufficient to consider the
analogous result for a configuration \(\config^n\) of \(n\)
points in the disk \(D_n\). For then we can apply the result to the
lesser configuration \(\by^{k(n)}\) (for \(k(n)\) as given in
condition (b) of Theorem \ref{thm:lower-bound}) and obtain
\[
\excess(G(\by^{k(n)}))\quad=\quad
\Omega(\sqrt{\log k(n)}) =
\Omega(\sqrt{ \log \pi\rho n \varepsilon})\quad=\quad\Omega(\sqrt{ \log n})\,,
\]
while
\begin{multline*}
\excess(G(\by^{k(n)}))\quad=\quad
  \frac{n(n-1)}{k(n)(k(n)-1)}\excess\left(G(\config^n)\right) \\
\quad\leq\quad
  \frac{1}{\pi\rho\varepsilon(\pi\rho\varepsilon-1/n)}
  \excess\left(G(\config^n)\right)\,,
\end{multline*}
from which Theorem \ref{thm:lower-bound} follows.

We therefore consider \(\config^n\subset D_n\) being \(L_n\)-equidistributed as the uniform distribution
on $D_n$.  So by definition there is a coupling \((X_1,Y_1)\)
(here we omit dependence on \(n\)) where \(X_1\) has uniform distribution on \(\config^n\),
\(Y_1\) has uniform distribution on \(D_n\) and
\begin{equation}
\Delta_n \quad=\quad
 \Expect{ \min\left(1, \frac{|X_1 - Y_1|}{L_n}\right)}
 \quad\to\quad 0 \mbox{ as } n \to \infty.
\end{equation}

Write \((X_2,Y_2)\) for an independent copy of  \(X_1,Y_1\).
In the definition of \emph{excess} it makes no asymptotic difference if we allow $j = i$ in
\(\average_{(i,j)}\), so we may take
\begin{equation}
  \excess(G(\config^n))\quad=\quad
  \Expect{\ell(X_1,X_2) - |X_1-X_2|} .
\end{equation}
Set
\begin{equation}\label{eq:equation17}
A_n \quad=\quad 
[|Y_1 - X_1| \leq L_n]\cap [|Y_2 - X_2| \leq L_n]
\end{equation}
so that by Markov's inequality
\begin{equation}
\label{eq:A_n}
\Prob{A_n} \quad\geq\quad 1 - 2 \Delta_n .
\end{equation}
Define \(\ell(Y_1, Y_2)\) by supposing that \(Y_i\) is plumbed in
to the network using a connection by a \emph{temporary} line segment with
endpoints \(Y_i\) and \(X_i\).
A direct computation  shows that on the event \(A_n\)
\begin{multline*}
  \ell(Y_1, Y_2)-|Y_1-Y_2| \quad\leq\quad\\
\left(\ell(X_1, X_2)+|X_1-Y_1|+|X_2-Y_2|\right)
-
\left(|X_1-X_2|-|X_1-Y_1|-|X_2-Y_2|\right)\\
\quad \leq \quad
\ell(X_1, X_2)-|X_1-X_2|
+4 L_n.
\end{multline*}
Consequently
\begin{equation}
  \label{eq:vasershtein-approximation}
  \Expect{\ell(Y_1, Y_2)-|Y_1-Y_2|; A_n} \quad\leq\quad
    \excess(G(\config^n)) + 4L_n\,.
\end{equation}
%
%
By hypothesis  \(L_n=o(\sqrt{\log n})\),
and so the proof of Theorem \ref{thm:lower-bound} reduces to showing
that the left side (the excess for two random points chosen uniformly
in the disk) is
\(\Omega(\sqrt{\log n})\).

\subsection{Near-parallelism for case of small excess}
\label{sec:lwb-geodesic-slope}

We now substantiate our previous remark that the network geodesic must
run almost parallel to the Euclidean geodesic if the excess is small.

It is convenient to situate the disk \(D_n\) in the complex plane
\(\mathbb{C}\) so as to have a compact notation for rotations.  For
\(t>0\) we define \(Z_t\) and \(\Phi\) by
\begin{align}
  \exp\left(i\Phi\right)&\quad=\quad
  \frac{Y_2-Y_1}{|Y_2-Y_1|}\,,\nonumber\\
  Z_t&\quad=\quad Y_1+
  t \times\exp\left(i\Phi\right)\,.
\label{eq:complex-construction}
\end{align}
%
%

\begin{figure}[ht]
  \centering
  \includegraphics[width=3in]{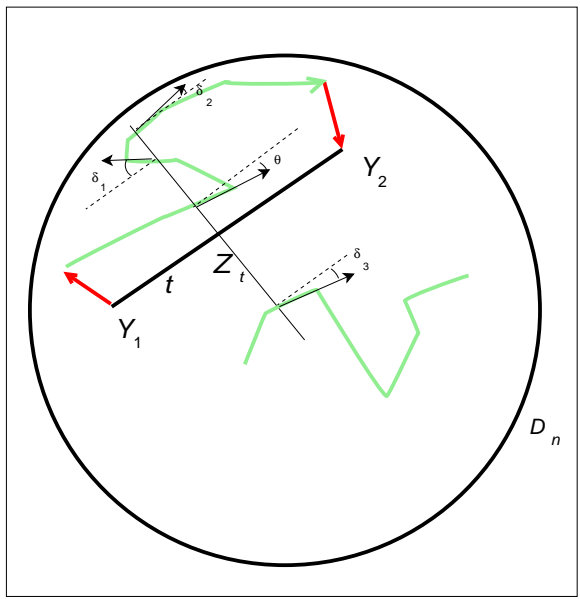}
  \caption{Illustration of construction of \(Y_1\), \(Y_2\), and
    \(Z_t\). The angles \(\theta(t)\) and \(\delta_1\), \(\delta_2\),
    \ldots\ are computed using the angles of incidence of network
    segments on the perpendicular running through \(Z_t\);
    \(\Upsilon_{t,\chi}\) is the minimum of absolute values of all
    such angles of points of intersection within
    \(\sqrt{2t\chi+\chi^2}\) of \(Z_t\).}
\label{fig:lwb-1}
\end{figure}

Let \(\gamma:[0,\ell(Y_1, Y_2)]\to\mathbb{C}\) be the unit-speed
network geodesic running from \(Y_1\) to \(Y_2\) (using the temporary
plumbing to move from \(Y_1\) to \(X_1\) and then again from \(Y_2\)
to \(X_2\)).  Then (bearing in mind that \(\left|\gamma^\prime(t)\right|=1\))
\begin{multline}
  \label{eq:slope-bound0}
  \ell(Y_1, Y_2)\quad=\quad
\int_0^{\ell(Y_1, Y_2)}\left|\gamma^\prime(s)\right|\d s
\quad\geq\quad
\int_0^{|Y_1-Y_2|}\left|\gamma^\prime(\tau(t))\right|\tau^\prime(t)\d t\,,
\end{multline}
where \(\tau(t)\) is the first time \(s\) at which
\(\langle\gamma(s)-Y_1,\exp\left(i\Phi\right)\rangle=t\). (Note that
our networks are formed from finite collections of line segments. Hence
\(\tau^\prime\) will be defined and finite save perhaps at a finite number of
times.) This and the following constructions are illustrated in
Figure \ref{fig:lwb-1}.

Defining \(\theta(t)\) by \(\sec\theta(t)=\tau^\prime(t)\), and
using \(\sec\theta\geq1+\tfrac{1}{2}\theta^2\), we deduce
\begin{equation}
  \label{eq:slope-bound}
  \ell(Y_1, Y_2)\quad\geq\quad |Y_1-Y_2|+
\half\int_0^{|Y_1-Y_2|}\theta(t)^2\d t\,.
\end{equation}

Furthermore we can use Pythagoras and the geodesic property of
Euclidean line segments to show the following. Let \(H(t)\) be the
maximum \(|r|\) for which, for some \(s\),
\[
\gamma(s) \quad=\quad Z_t + i r \exp\left(i\Phi\right)\,.
\]
If the excess for the network
geodesic from \(Y_1\) to \(Y_2\) is bounded above by \(\ell(Y_1,
Y_2)-|Y_1-Y_2|\leq \chi\) then \(H(t)\leq\sqrt{2t\chi+\chi^2}\).

Let \(\Upsilon_{t,\chi}\) be the smallest \(|\delta|\) such that some network
segment intersects the perpendicular \(\{Z_t+ir\exp{i\Phi}:r\in\Reals\}\) at angle \(\pi/2+\delta\) and at
distance at most \(\sqrt{2t\chi+\chi^2}\) from
\(Z_t\) (thus \(\delta\) is the angle of incidence of this network
segment on the perpendicular).
If
\(\ell(Y_1, Y_2)-|Y_1-Y_2|\leq \chi\) and
\(|Y_1-Y_2|\geq\kappa\sqrt{\rho n}\), we can use (\ref{eq:slope-bound}) to deduce
\begin{multline*}
  \ell(Y_1, Y_2)-|Y_1-Y_2|\quad\geq\quad\\
\half\int_0^{\kappa\sqrt{\rho n}}\Upsilon_{t,\chi}^2\d t
-\half\left(\frac{\pi^2}{4}\right)\times\left(|X_1-Y_1|+|X_2-Y_2|\right)
\,.
\end{multline*}
(The second summand allows for the temporary plumbing in of
connections \(X_1Y_1\) and \(X_2Y_2\), for which the angle
\(\theta(t)\in(0,\tfrac{\pi}{2})\) is not controlled by permanent
network segments).
So introduce the event
\begin{equation}
  B_{\kappa, \chi}\quad=\quad\left[
\ell(Y_1, Y_2)-|Y_1-Y_2|\leq \chi\,,
|Y_1-Y_2|\geq\kappa\sqrt{\rho n}
\right]\,
\label{eq:restriction}
\end{equation}
and recall from Equation (\ref{eq:equation17}) the event
\(A_n = \cap_{i=1}^2 [|Y_i - X_i| \leq L_n]\).
Taking expectations, we deduce
\begin{multline*}
  \Expect{\ell(Y_1, Y_2)-|Y_1-Y_2|\;;\;B_{\kappa, \chi}\cap A_n}\\
  \quad\geq\quad
\half\int_0^{\kappa\sqrt{\rho n}}
\Expect{\Upsilon_{t,\chi}^2 \;;\;B_{\kappa, \chi}\cap A_n
}\d t
-\frac{\pi^2}{4}L_n\,.
\end{multline*}
Using integration by parts to replace the expectation by a
probability,
\begin{multline}
    \label{eq:slope-lower-bound}
  \Expect{\ell(Y_1, Y_2)-|Y_1-Y_2|\;;\;B_{\kappa,
      \chi}\cap A_n}+\frac{\pi^2}{4}L_n\\
\quad\geq\quad
\int_0^{\kappa\sqrt{\rho n}}\int_0^\infty
\Prob{[\Upsilon_{t,\chi}>u]\cap B_{\kappa, \chi}\cap A_n
}u\,\d u\d t\\
\quad = \quad
\int_0^{\kappa\sqrt{\rho n}}\int_0^\infty
\left(\Prob{B_{\kappa, \chi}\cap A_n} - \Prob{[\Upsilon_{t,\chi}\leq u]\cap B_{\kappa, \chi}\cap A_n
}\right)u\,\d u\d t \\
\quad \geq \quad
\int_0^{\kappa\sqrt{\rho n}}\int_0^\infty
\max \left(\Prob{B_{\kappa, \chi}\cap A_n} - \Prob{\Upsilon_{t,\chi}\leq u
}, 0\right)u\,\d u\d t\, .
\end{multline}
Note also that from the definitions of \(B_{\kappa, \chi}\) and \(A_n\), using
\eqref{eq:A_n}, \eqref{eq:vasershtein-approximation} and Markov's inequality
\begin{multline}
\label{eq:BA}
1 - \Prob{B_{\kappa, \chi}\cap A_n} \quad=\quad
1-\Prob{A_n} + \Prob{A_n\setminus B_{\kappa, \chi}}
\\
\quad\leq\quad
 2\Delta_n + \Prob{|Y_1-Y_2|<\kappa\sqrt{\rho n}}
+ \frac{ \excess(G(\config^n)) + 4L_n}{\chi} .
\end{multline}
To make progress we now need to find an upper bound for
\(
\Prob{\Upsilon_{t,\chi}\leq u}\)
and this is the subject of the next section.

\subsection{Upper bounds using uniform random variables }
\label{sec:lwb-uniform}

\begin{figure}[ht]
\centering
    \includegraphics[width=3in]{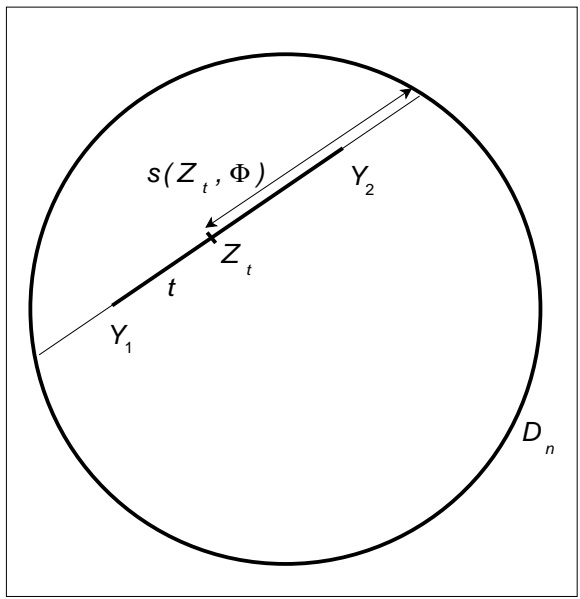}
    \caption{Illustration of construction in Lemma \ref{lem:polar-control}.}
\label{fig:lwb-2}
  \end{figure}

Firstly we compute an upper bound on the joint density of the
quantities \({Z}_t\) and \({\Phi}\) from the previous section,
illustrated in Figure \ref{fig:lwb-2}.

  \begin{lem}\label{lem:polar-control}
    Suppose \(Y_1\), \(Y_2\) are independent
    uniformly distributed random points in a disk \(D_n\) of radius
    \(\sqrt{\rho n}\) and centre \(0\) in the complex plane
    \(\mathbb{C}\). With \(Z_t\) and
    \(\Phi\) defined as in \eqref{eq:complex-construction}, the
     joint density of \(Z_t\) and \(\Phi\)
    is given over \(\mathbb{C}\times[0,2\pi)\) by
    \begin{equation}
      \Indicator{z-t e_\phi\in D_n} \frac{\left(t+s(z,\phi)\right)^2}{2\pi^2\rho^2n^2}
      {\d}z\,\d \phi\,,
      \label{eq:joint-density1}
    \end{equation}
    where \(e_\phi=e^{i\phi}\) is the unit vector making angle
    \(\phi\) with a reference \(x\)-axis, and \(s(z,\phi)\)
is the distance from \(z\) to
the disk boundary \(\partial D_n\) in the direction \(\phi\) (thus in particular
\(z+s(z,\phi)e_\phi\) is on the disk boundary).
\end{lem}
\begin{proof}
  Express the joint density for \(Y_1\),
  \(Y_2\) as a product of a uniform density over \(D_n\) for
  \(Y_1\) and polar coordinates \(r,\phi\) about \(Y_1\)
  for \(Y_2\):
\[
\Indicator{y_1\in D_n} \frac{{\d}y_1}{\pi\rho n}
\Indicator{y_1+r e^{i\phi}\in D_n} \frac{r\,\d r\,\d \phi}{\pi\rho n}\,.
\]
Obtain the result by integrating out the \(r\) variable and
transforming the \(y_1\) variable to \(z\) by \(z=y_1+t e^{i\phi}\).
\end{proof}
\begin{corollary}\label{lem:polar-control-corollary}
  The density for \(Z_t\) and
  \(\Phi\pmod\pi\) is
  \begin{multline}
    \label{eq:joint-density}
f(z,\phi)\quad=\quad\\
\left(
\Indicator{z-t e_\phi\in D_n} \frac{\left(t+s(z,\phi)\right)^2}{2}
+
\Indicator{z+t e_\phi\in D_n} \frac{\left(t+s(z,\pi+\phi)\right)^2}{2}
\right)\times\\
 \times\Indicator{0\leq\phi<\pi}
 \frac{{\d}z\,\d \phi}{\pi^2\rho^2n^2}\,.
  \end{multline}
with an upper bound
\begin{equation}
  \label{eq:polar-control-upper-bound}
f(z,\phi)\quad\leq\quad
4 \times \Indicator{0\leq\phi<\pi}
\frac{{\d}z\,{\d}\phi}{\pi^2\rho n}\,.
\end{equation}
\end{corollary}
\begin{proof}
  Equation \eqref{eq:joint-density} follows immediately from adding
  the two expressions from Equation \eqref{eq:joint-density1} for
  \(\phi\pmod\pi\). The upper bound follows by noting
  \begin{enumerate}
  \item the maximum will occur when \(z-t e_\phi\) runs along a
    diameter as \(t\) varies;
\item furthermore when one of \(z\pm t e_\phi\) lies on the disk
  boundary;
\item and furthermore when \(z=0\) is located at the centre of the disk
  (so \(t=s(z,\pm\phi)=\sqrt{\rho n}\)).
  \end{enumerate}
\end{proof}

Now consider the line segment \({\mathscr{S}}_{t,\chi}\) centred at
\(Z_t\), with end-points given by the pair \(\pm
i\sqrt{2t\chi+\chi^2}\exp\left(i\Phi\right)\);
and consider the rose-of-directions empirical measure of
 angles made by
intersections of network edges with this segment:
\begin{equation}
  \label{eq:rose-of-directions}
{\mathscr{R}}_{t,\chi}(A)\quad=\quad
\Number\left\{\text{ network intersections on
  }{\mathscr{S}}_{t,\chi}
\text{ with angle of incidence lying in } A\right\}
\end{equation}
(here angles are measured modulo \(\pi\), and \(A\subseteq[0,\pi)\)).
%
%
We may apply a Buffon-type argument to bound
\(\Expect{{\mathscr{R}}_{t,\chi}(A)}\) using Inequality
\eqref{eq:polar-control-upper-bound}. Consider the contribution to the
expectation from a fixed line segment of the network of length
\(\ell\): the result of disintegrating the integral expression for
this according to the value of \(\phi\) is an integral of
\(f(z,\phi)\) with respect to \(z\) over a region formed by
intersecting the disk with a parallelogram of base side-length
\(\ell\) and height \(2\sqrt{2t\chi+\chi^2}\sin\alpha\) (here
the angle \(\alpha\) depends implicitly on \(\phi\) and \(z\)). Of
course the integral vanishes if \(\phi\not\in A\). Thus Inequality
\eqref{eq:polar-control-upper-bound} yields a bound
\[
\Expect{{\mathscr{R}}_{t,\chi}(A)}\quad\leq\quad
\frac{4 }{\pi^2\rho n}\times
\int_{G(\config^n)}\int_A 2 \sqrt{2t\chi+\chi^2} \sin\alpha
\d \alpha \d z.
\]
For constant  \(\chi\), the event \( [{\Upsilon}_{t,\chi}\leq u ]\) is the event
\([{\mathscr{R}}_{t,\chi}(\tfrac{\pi}{2} - u, \tfrac{\pi}{2} + u) \geq 1]\) and so
\begin{equation}
\label{eq:uniform-estimate}
  \Prob{{\Upsilon}_{t,\chi}\leq u}
\quad\leq\quad
\Expect{{\mathscr{R}}_{t,\chi}(\tfrac{\pi}{2} -u, \tfrac{\pi}{2} + u)}
\quad\leq\quad
\frac{16}{\pi^2\rho}\frac{\len(G(\config^n))}{n}\sqrt{2t\chi+\chi^2} \times u\,.
\end{equation}

\subsection{Calculations}
\label{sec:calculations}
We have assembled the ingredients for the proof of Theorem
\ref{thm:lower-bound}, and so now 
can perform the calculations to get a quantitative lower bound.

We proceed by contradiction. Suppose that \(\excess(G(\config^n))= o(\sqrt{\log n})\).
Inspecting \eqref{eq:BA} we see that we can choose \(\chi=\chi_n = o(\sqrt{\log n})\)
and some small $\kappa > 0$ such that
 for all sufficiently large $n$
\begin{equation}
 \label{kappa-chi}
 \Prob{B_{\kappa, \chi}\cap A_n} \geq 2^{-1/3} .
\end{equation}
So we can combine \eqref{eq:vasershtein-approximation}  and
\eqref{eq:slope-lower-bound}
(and the fact that \(\pi^2/4<3\)) to get
\[  \excess(G(\config^n)) + 7L_n \geq
\int_0^{\kappa\sqrt{\rho n}}\int_0^\infty
\max \left(2^{-1/3} - \Prob{\Upsilon_{t,\chi}\leq u
}, 0\right)u\,\d u\,\d t . \]
By \eqref{eq:uniform-estimate} and hypothesis of Theorem \ref{thm:lower-bound}, there exists a constant $B$ such that
\[   \Prob{{\Upsilon}_{t,\chi}\leq u}
\quad\leq\quad \sqrt{\frac{B}{12}} \  \sqrt{2t\chi+\chi^2} \times u . \]
Applying the formula
\( \int_0^\infty \max(0, \alpha - \beta u) u \,\d u = \frac{\alpha^3}{6\beta^2}\)
we see
\begin{equation}
  \excess(G(\config^n)) + 7L_n \geq \frac{1}{B}
\int_0^{\kappa\sqrt{\rho n}} \frac{1}{2t\chi+\chi^2}\, \d t
= \frac{\log(\kappa\sqrt{\rho n} + \tfrac{\chi}{2}) - \log \tfrac{\chi}{2} }{2\chi B}.
\label{eq:Gx}
\end{equation}
Recall this holds under the assumption that $\chi_n= o(\sqrt{\log n})$
and that $\kappa>0$ is constant.
We are given that
\(L_n = o(\sqrt{\log n})\), and we have supposed for the purposes of
contradiction that \(\excess(G(\config^n))= o(\sqrt{\log n})\).
But then \eqref{eq:Gx} takes the form
\[
 o(\sqrt{\log n}) \quad\geq\quad \frac{\Omega(\log n)}{o(\sqrt{\log n})}\,,
\]
which is impossible. We deduce we must have \(\excess(G(\config^n))= \Omega(\sqrt{\log n})\).

 \section{Closing remarks and supplements}\label{sec:conclusion}
\subsection{Spatial network design}\label{sec:spatial-network-design}
 Within the realm of spatial network design, the closest work we know
 is that of Gastner and Newman \cite{GastnerNewman-2006}, who consider the
 similar notion of a {\em distribution network} for transporting
 material from one central vertex to all other vertices.  They give a
 simulation study (their Figure \(2\)) of a certain algorithm on random
 points, and comment
 \begin{quote}
   Thus, it appears to be possible to grow networks that cost only a
   little more than the [minimum-length] network, but which have far
   less circuitous routes.
 \end{quote}
 Our Theorem \ref{thm:upper-bound} provides a strong formalisation of this idea.

 An algorithm for minimizing excess for a given length is described in
 \cite{SchweitzerEbelingRoseWeiss-1998}, where results for a 39 point
 configuration are shown.  But neither this nor
 \cite{GastnerNewman-2006} has led to study of $n \to \infty$
 asymptotics.

\subsection{Fractal structure of the Steiner tree on random points}\label{sec:fractal}
A longstanding topic of interest in statistical physics is that of
 the continuum limits of
various discrete two-dimensional self-avoiding walks arising in
probability models, \emph{eg}
\begin{itemize}
\item uniform self-avoiding walks on the lattice,
\item paths within uniform spanning trees in the lattice,
\item paths within minimum spanning trees in the lattice.
\end{itemize}
This study
has recently been complemented by spectacular successes of rigorous
theory \cite{LawlerSchrammWerner-2004}.  It is conjectured that routes
in Steiner trees on random points have similar fractal properties
\cite{Read-2005}: route-length between points at distance \(n\)
should grow as \(n^\gamma\) for some \(\gamma > 1\).  However, as our
construction shows, such results have little relevance to spatial
network design.

\subsection{The counterintuitive observation}\label{sec:counterintuitive}
The counterintuitive observation following Definition
\ref{def:ratio-statistic} follows quickly from the work of Theorem
\ref{thm:upper-bound}. Suppose the configuration \(\config^n\) is
well-dispersed, in the weak sense that for some \(\gamma\in(0,1)\) we
find the number of point pairs within \(n^{\gamma/2}\) of each
other is \(o\left(\tbinom{n}{2} n^{\gamma-1}\right)\) (certainly this is
the case for most patterns generated by uniform random sampling from
\([0,\sqrt{n}]^2\)). Consider a network \(G(\config^n)\) produced by
augmenting the Steiner tree according to the construction in the proof
of Theorem \ref{thm:upper-bound}. Using the properties of this
construction, the following can be shown
\begin{align*}
  \Expect{\ratio\left(G(\config^n)\right)} &\quad=\quad
\Expect{\average_{(i,j)}\frac{\ell(x_i,x_j)}{|x_i-x_j|}-1}\\
&\quad\leq\quad \text{constant } \times o(n^{\gamma-1}) +
(1-o(n^{\gamma-1})) \left(\frac{O(\log
    \sqrt{2n})}{n^{\gamma/2}}\right)\\
&\quad\leq\quad
O\left(\max\left(\frac{1}{n^{1-\gamma}},\frac{\log n}{n^{\gamma/2}}\right)\right)\,.
\end{align*}

\subsection{Derandomisation}\label{derandomization}
Theorem \ref{thm:upper-bound} is a purely deterministic assertion, though
our proof used randomisation  (supplied by the Poisson line
process).  It seems intuitively plausible that one could give a purely
deterministic proof, say by replacing the Poisson line process with a
suitable sparse set of deterministically positioned lines having a
dense set of orientations.

\subsection{Quantifying equidistribution}\label{sec:equidistribution}
 The classical equidistribution property
 \begin{quote}
   the empirical distribution of \(\{n^{-1/2}x^n_i, 1 \leq i \leq
   n\}\) converges weakly to the uniform distribution on \([0,1]^2\)
 \end{quote}
 is equivalent (by a straightforward argument) to the property
 \begin{quote}
 \(\config^n\) is \(L_n\)-equidistributed as the uniform distribution on the square of area $n$,
 for \emph{some} \(L_n = o(n^{1/2})\).
 \end{quote}
  Replacing one sequence
 of \(L_n\) by another slower-growing sequence makes
equidistribution a stronger assumption, and so our
 assumption in Theorem \ref{thm:lower-bound}(a)  (equidistribution
 for some \(L_n = o(\log^{1/2}
 n)\)) is stronger than the classical equidistribution property.
 Indeed Theorem \ref{thm:lower-bound} fails under the classical
 equidistribution property, as the following example shows.
 \begin{ex}
   Let \(L_n = n^\gamma\) for some  \(\gamma\in\left(\tfrac{3}{8}, \tfrac{1}{2}\right)\).  There
   exist networks \(G(\config^n)\) which are
   \(L_n\)-equidistributed as the uniform distribution on the square of area $n$,
 for which \(\len(G(\config^n)) = o(n)\) whilst
   \(\excess(G(\config^n))\to0\).
\end{ex}
For example: partition \([0,n^{1/2}]^2\) into subsquares of side
\(L_n/\log n\), construct the complete graph on all centres of such subsquares,
allocate the \(n\) points evenly amongst subsquares and position them
arbitrarily close to the centres.

As is apparent from the non-stochastic condition implying
\(L_n\)-equidistribution, there is a wide variety of configurations
satisfying \(L_n\)-equidistribution. Here we consider the particular
case of independent uniform sampling, and show that this generates an
\(L_n\)-equidistributed sequence of configurations.
\begin{remark}\label{rem:reassuring-exercise}
  Sample the configuration \(\config^n\) independently and uniformly
  from \([0,\sqrt{n}]^2\).  Let \(L_n \to \infty\), perhaps arbitrarily
  slowly.  Then the probability that the configuration \(\config^n\)
  is \(L_n\)-equidistributed with the uniform distribution converges
  to \(1\). This follows by dividing \([0,\sqrt{n}]^2\) into cells of
  side-length asymptotic to \(L_n/\sqrt{2}\), by conditioning on
  \(\config^n\), and by ``blurring'' the points of \(\config^n\) by
  replacing each point \(x\in\config^n\) by an independent draw taken
  uniformly from the cell containing \(x\). Then a uniform random draw
  \(\widetilde{Y}_n\) of one of the blurred points can be coupled to
  lie within \(L_n\) of a uniform random draw \(X_n\) from the
finite configuration
  \(\config^n\). A simple argument using the Binomial distribution
  then shows that the total variation distance between
  \(\widetilde{Y}_n\) and Uniform(\([0,\sqrt{n}]^2\)) tends to zero;
  it follows that \(X_n\) can be coupled to a
  Uniform(\([0,\sqrt{n}]^2\)) random variable \(Y_n\) so that
\[
\Expect{\min\left(1, \frac{|X_n-Y_n|}{L_n}\right)|\config^n}\quad\to\quad0\,,
\]
where the convergence takes place in probability.
%
%
\end{remark}

\subsection{Poisson line process networks}\label{sec:poisson-line-networks}
Remark \ref{rem:further-asymptotic-mean-cell-boundary} indicates that
more can be said about the mean semi-perimeter
\[
\tfrac{1}{2}\Expect{\len(\partial\mathcal{C}(v_i,v_j))}\,,
\]
 and this
will be returned to in later work. For example, consider the network
formed entirely from a Poisson line pattern. If the pattern is
conditioned to contain points \(v_i\), \(v_j\) then the perimeter
\(\partial\mathcal{C}(v_i,v_j)\) will be close to providing a genuine network
geodesic.

Note that questions about \(\mathcal{C}(v_i,v_j)\) bear a family resemblance to the D.G.~Kendall conjecture about the
asymptotic shape of large cells in a Poisson line pattern. However
\(\mathcal{C}(v_i,v_j)\) is the result of a very explicit
conditioning and hence explicit and rather complete
answers can be obtained by direct methods, in contrast to the striking
work on resolving the conjecture about large cells
\cite{Miles-1995,Kovalenko-1997,Kovalenko-1999,HugReitznerSchneider-2003}.

\subsection{An open question}
In the random points model we can pose a more precise question.  Over
choices of network \(G\) subject to the constraint
\[
\Expect{\len(G(\config^n)) - \len(\Steiner(\config^n))} \quad=\quad o(n)\,,
\]
 or the constraint
\[
\Expect{\len(G(\config^n)} \quad=\quad O(n)\,,
\]
 what is the minimum value of \(\Expect{\excess(G(\config^n))}\)?
Our results pin down this minimum value, in the latter case to the range
\([\Omega\left(\sqrt{\log n}\right), O(\log n)]\)
   and in the former case the range
\([\Omega\left(\sqrt{\log n}\right), o(w_n \log n)]\).
  But it remains an open question what should be the exact order of magnitude.

\bigskip

\begin{tabbing}
David J. Aldous, \=
Department of Statistics, \=
367 Evans Hall \#\  3860, \\
\>U.C. Berkeley CA 94720, USA\\
Email: \>\url{aldous@stat.berkeley.edu}\\
URL:   \>\url{www.stat.berkeley.edu/users/aldous}\\
\\
Wilfrid Kendall, \>
Department of Statistics, \\
\>University of Warwick, Coventry CV4 7AL, UK\\
Email: \>\url{w.s.kendall@warwick.ac.uk}\\
URL:   \>\url{http://www.warwick.ac.uk/go/wsk}
\end{tabbing}


   \bibliographystyle{plain}
   \bibliography{networks}
   \label{sec:bibliography}


\end{document}